\documentclass[letterpaper, 10 pt, conference]{ieeeconf}  
\usepackage[hidelinks]{hyperref}
\usepackage{amsmath, amssymb, amsthm, latexsym, graphics, graphpap, layout, wasysym, multicol, enumerate,amsfonts,mathrsfs}
\usepackage[pdftex]{graphicx}
\usepackage{caption}
\usepackage{xcolor}
\IEEEoverridecommandlockouts                              
\title{\LARGE {\bf{
CAsimulations: Modeling of topological dynamics in a disease using cellular automata}}}

\author{Jorge Andrés Ibáñez Huertas$^1$, Carlos Isaac Zainea Maya$^2$ \\
$^1$Departamento de Matemáticas, Universidad Central, Bogotá, Colombia. jibanezh@ucentral.edu.co\\
$^2$Departamento de Matemáticas, Universidad Externado de Colombia, Bogotá, Colombia.\\ carlos.zainea@uexternado.edu.co}

\begin{document}

\maketitle
\thispagestyle{empty}
\pagestyle{empty}

\begin{abstract}


The prediction of the behavior of the disease, the level of affectation in a population and the ways to control it are the most important aspects studied by epidemiology using tools such as historical data and mathematical models.
So, our objective is (1) to provide a tool capable of analyzing epidemiological phenomena starting from the most common social interactions within a group of individuals.
(2) To provide a methodology to build epidemiological models from patterns and logical rules. 
(3) Determine the impact of social interactions on the spread of a disease.
This paper describes the logical construction of two epidemiological models in cellular automata together with two of their variations based on topological and dynamical principles.

\end{abstract}

\section{INTRODUCTION}

Today, there are many mathematical models to project the progress of a disease. Contingencies such as the one experienced with COVID-19 promote a greater understanding of these models, and new theories are created around this problem from different perspectives. Among these models, the compartmental models allow studying the evolution of the different states considered. For example, with SIR the evolution of susceptible, infected and recovered individuals is studied separately; A peculiarity of this model is that it projects the growth of individuals in each of its compartments, and it is possible to establish at what point we are before the disease becomes endemic, however, it does not take into account the possible factors related to the care, risk or exposure of some individuals.

In this work we will show, based on simulation techniques, a possible approach to compartmental models that take into account the closeness with other individuals. To achieve this, we use cellular automata theory and an extension of it, using topological neighborhood concepts and fundamental neighborhood systems. In {\bf{theorem} [1]} we show how the impact of one cell on the others constitutes a fundamental system of neighborhoods, and after we define some evolution rules for SI, SIS and SIR models for each cell. Finally, an example that allows us to study the models described in this paper in a village of 49 people.

In the development of this work, a Python package called CAsimulations was made, the documentation of this package can be consulted at \href{https://grupo-de-simulacion-con-automatas.github.io/CAsimulations-Modelacion-de-dinamicas-topologicas-en-la-propagacion-de-una-enfermedad-usando-CA/}{\color{blue} https://grupo-de-simulacion-con-automatas.github.io/CAsimulations-Modelacion-de-dinamicas-topologicas-en-la-propagacion-de-una-enfermedad-usando-CA/} 

\section{PRELIMINARIES}
\subsection{EPIDEMIOLOGICAL STUDY}
One of the most important objects of study in epidemiology is to determine the endemic nature of the disease, that is, whether the disease will affect the population for a long time or whether it will gradually disappear. Usually, this is determined by the indicator $\mathcal{R}_0$ known as the basic reproduction number and which represents the number of individuals that were infected by patient zero in a susceptible population. Normally, if $\mathcal{R}_0<1$ the disease will slowly disappear and if $\mathcal{R}_0>1$, we may have an endemic case.

Heesterbeek and Dietz defining the basic reproduction number as:

\begin{equation}\label{eq:R0}
 \mathcal{R}_0 = \int_0^\infty b(t)F(t) dt,
\end{equation}

where $b(t)$ represents the average number of new infections that an infected individual will produce during time $t$ and $F(t)$, known as the survival function, represents the probability that a newly infected individual will remain in that state for at least time $t$ \cite{conceptOfR0, perspectivesOnR0}.

Generally, when we talk about epidemiological models, we consider three states or classes into which we can divide the population over time: those who can contract the infection, those who are infected and those who are recovered. If the recovered populations do not acquire permanent immunity, we have a SIS model, in which the susceptible population can contract the infection and once they recover they come back to susceptible state. On another way, if the recovered individuals generate disease immunity, we have an SIR model \cite{mateModelsInPopulationAndEpidemiology}.

\subsection{THE \textit{SIS} MODEL}

The SIS model considers two possible states, susceptible (S) and infected (I). Variations between states are given by new infections and individuals recovering from the disease. In addition, each state is affected by the mortality/birth rate and/or the parameter describing death caused by the disease.

Normally, when talking about epidemiological models with disease deaths, four parameters are considered: 

\begin{itemize}
    \item The \textbf{infection rate $\beta$}, which represents the probability of a susceptible individual acquiring the disease after contact with an infected individual.
    \item The \textbf{recovery rate $\alpha$}, which can be understood as the probability that an infected person will recover from the disease. It sometimes represents the average time it takes for an infected person to recover from the disease \cite{diego2010}.
    \item The \textbf{birth/mortality rate $\mu$}, which in the case of traditional models are assumed to be equal. The birth rate indicates the number of individuals entering the space, and the mortality rate represents the individuals who die from causes other than disease.
    \item The \textbf{death rate due to disease $\theta$}, which indicates the probability that an infected person is likely to die from the disease.
\end{itemize}

In our work, we will use a constant population and thus $S'+I'=0$. With this in mind, we can describe the model from a system of differential equations as follows:

\begin{equation}\label{eq:modeloSIS}
\left\{
\begin{array}{l}
S' = \mu(1 - S) + (1 - \theta)\alpha I - \beta S I, \\
I' = \beta S I - (1 - \theta)\alpha I - \mu I.
\end{array}
\right.
\end{equation}

In \cite{diego2010} they use the equation \ref{eq:R0} to determine in an analogous way that $\mathcal{R}_0=\frac{\beta}{\alpha(1-\theta)+\mu}$ so we will not go deeper into this process.

\subsection{THE \textit{SIR} MODEL}

For this model, the state of immunity (R) to the disease is considered. Unlike the SIS model, in the SIR model there is no interaction from state I to state S, since it is assumed that individuals who recover from the disease will not be able to contract it again, so they will pass to state R. In this way,

\begin{equation}\label{eq:Modelo SIR}
\left\{
\begin{array}{l}
S' = \mu(1 - S) + \alpha\theta I - \beta S I, \\
I' = \beta S I - \alpha I - \mu I\text{, and } \\
R' = \alpha I - \alpha\theta I - \mu R.
\end{array}
\right.
\end{equation}

Using an argument similar to the one used to determine the value of $\mathcal{R}_0$ in the SIS model, it is possible to determine that for the SIR model, $\mathcal{R}_0 = \frac{\beta}{\alpha+\mu}.$ 

\subsection{CELLULAR AUTOMATAS}

Cellular automata (or CA) were born with the work of Von Neumann in the late 1940s with his paper "The General and Logical Theory of Automata" in which he first put forward ideas for a machine capable of self-replication. He worked on a discrete two-dimensional system to develop fairly complex dynamics that were also self-replicating. \cite{alfons2010,ACaplicacionesComputacion}.

We will now discuss some general information about cellular automata that will be useful in the following sections.

\begin{itemize}
    \item \textbf{Cell space $\mathcal{L}$:} Is the set where all the cells considered for the model live and interact. In general this space is discrete, regular and finite, the latter due to the computational limitations present in the tools with which models are built in cellular automata.
    \item \textbf{Set of states $\Sigma$:} Is the finite set of all possible categories in which the cells of the space $\mathcal{L}$ can be. Each element $\sigma$ of $\Sigma$ will be known as a state of the model.
    \item \textbf{Neighborhoods:} It is the set of cells "close" to a particular one. In general, we do not work with the whole set of neighborhoods $\mathcal{V}(x)$ but instead, we consider elements of each of these families to form a set of neighborhoods over the space $\mathcal{L}$. This set is known as a \textit{neighborhood system} over $\mathcal{L}$.
    \item \textbf{Evolution rules:} Evolution rules define the way in which the states of each cell change, taking into account the state of its neighbors. In general, we define an $\phi$ evolution rule as follows:
    $$\phi:\Sigma_x\times\overbrace{\Sigma\times\Sigma\times\cdots\times\Sigma}^{N}\longrightarrow\Sigma_x,$$
    where $N$ is the number of neighbors of $x$ and $\Sigma_x$ the set of states that $x$ can take. 
    To determine the evolution of the space $\mathcal{L}$ one must apply the evolution rule simultaneously on each of its cells.
\end{itemize}

It is possible to identify certain types of neighborhoods for each of the cells. For the purposes of this paper we will focus only on the minimal neighborhoods and some of its properties, for this we will use concepts set forth in \cite{elementosTopologiaGeneral, barmak2011, munkres, NeiraNacional}.

\textbf{Definition:} A minimal neighborhood of $x$ is the intersection of all the neighborhoods of $x$. We denote it as $U_x$.

\section{EPIDEMIOLOGICAL MODELS IN CA}

Let us think for a moment that if an individual susceptible to a disease has contact with many infected individuals, it can get sick much more easily than an individual who has contact with few infected individuals. Now, how many interactions with infected individuals, out of all those that a cell can have, are enough to generate a considerable probability of contagion? To answer this we must first determine the nature behind the possible relationships between cells, and for this we will make use of tools such as the fundamental systems of neighborhoods.

\subsection{SOCIAL INTERACTIONS AND IMPACTS}

Unlike the work done in \cite{populationDensity} in which each cell represented a region, we will consider each division as a single individual that will be provided with a set of qualities such as health status, age, neighbors, etc. These qualities will be given by the needs of the model we are developing, for example in the simplest SIR and SIS models it will not be necessary to endow the cells with ages, but we will have to take into account their neighborhoods and their health states.

Let us now think about the characteristics behind the close relationship between individuals or cells. We will denote the relationship between cells with the symbol, $\thicksim$ and having said that we have that:

\begin{itemize}
    \item All cells are in contact with themselves, so for each cell $x$ it is satisfied $x \thicksim x$.
    \item If a cell were in contact with some other cell, then that cell would be in contact with the first cell, i.e., $x\thicksim y$ implies $y\thicksim x$.
    \item If a cell interacts with two other cells, it does not necessarily imply that they interact with each other, so $x\thicksim y$ and $x\thicksim z$ do not imply that $y\thicksim z$.
\end{itemize}

Let us consider for example the set $X=\{a,b,c\}$ and the topology $\mathcal{A}=\{\emptyset,\{a\},\{b\},\{a,b\},\{a,c\},\{a,b,c\}\}$. The interaction relations present in $\mathcal{A}$ are:
$$\begin{array}{cccccc}
    a\thicksim a, & b\thicksim b, & a\thicksim b, & a\thicksim c & \text{, and} & b\thicksim c,
\end{array}$$
together with their equivalent symmetric relations.

The notion of relation in the previous example of cell interaction can be taken a bit further. Let us think for a moment about the impact that a behavior on cell $b$ can have on cell $a$, although these two cells do not interact with each other, the relationship that each of them has with cells $c$ and $d$ can have an impact on the other. Whereupon, we define the following interaction relationship:

\textbf{Definition:} We define the degree of impact between two points $a$ and $b$ as the least amount of interactions needed to get from $a$ to $b$. To recognize the degree of impact between two points, we will use the notation $a\thicksim_n b$ where $n\in\mathbb{N}$ denotes the least amount of interactions between $a$ and $b$.

If we return to the example, we can identify the following degrees of impact:

$$\begin{array}{cccc}
    a\thicksim_0 a, & b\thicksim_0 b, & a\thicksim_1 c, & a\thicksim_1 b, \\
    b\thicksim_1 a, & b\thicksim_2 c, & c\thicksim_0 a
\end{array}$$

\textbf{Theorem: \label{theorem}} The degrees of impact of an $x$ cell define a fundamental system of neighborhoods.

\begin{proof}
Let $x\in\mathcal{L}$ and $\mathcal{V}(x)$ be a family of neighborhoods of $x$. Define the set $A_0$ as the set of points with degree of impact with $x$ is equal to zero and recursively to the sets $A_k$, whose elements have degree of impact with $x$ is equal to or less than $k$. Clearly, $A_i\subseteq A_j$ for $0\leq i\leq j$ and thus $A_i\in\mathcal{V}(x)$ for $i=0,1,\cdots,n$.

Let us now see that the set $A_0$ is the minimal neighborhood of $x$. Consider $y\in U_x$, in particular, $y\in\bigcap\mathcal{V}(x)$. If the degree of impact of $y$ with $x$ is greater than zero by definition of impact grade, we can assert that there exists $z\in\mathcal{L}$ such that $z\thicksim x$, $z\thicksim y$ and $x\not\thicksim y$ and thus $x$ and $y$ are separable points. 

Consider $y\in A_0$ and by the definition of zero impact degree, we claim that the points $x$ and $y$ are not separable, so $y\in \mathcal{V}(x)$ and with this we conclude the proof.$\Box$
\end{proof}

One thing to keep in mind is that the degree of impact alone does not provide us with a measure of the impact of state changes of "distant" cells (or of the degree of impact greater than zero). We can define and understand these impact measures as the probability that a state change will affect the cell with which we are making the comparison. Thus, the impact indices will be values between 0 and 1 that can be given by any type of function that has as its domain the set of impact degrees.

\subsection{EVOLUTION RULES}

Before we start defining our evolution rules, we will define the following notations:
\begin{itemize}
    \item The state of the cell $x$ at a time $t$ will be denoted as $\pi^t(x)$.
    \item The number of individuals with impact degree $g$ and state $K$ of a cell $x$ at a time $t$, will be represented as $\sigma_{g,K}^t(x)$.
    \item To represent the number of individuals with an impact degree $g$ we will use the symbol $\Delta_g$. 
    \item We will use the symbols $\mathcal{S}^t,\mathcal{I}^t,\mathcal{R}^t$ and $\mathcal{D}^t$ to denote the sets of susceptible, infected, recovered and dead cells respectively in the space $\mathcal{L}$ at time $t$. Formally
    $$\mathcal{S}^t=\{x\in\mathcal{L}:\pi^t(x)=S\},$$
    and analogously, we define the sets $\mathcal{I}^t,\mathcal{R}^t$ and $\mathcal{D}^t$. Note that $$\mathcal{S}^t\cup\mathcal{I}^t\cup\mathcal{R}^t\cup\mathcal{D}^t=\mathcal{L}\text{ for all time }t.$$
\end{itemize}

In the following, we will define the rules for the simplest versions of the SIS and SIR models, for this we have to take into account that both models share the passage from state S to state I so we will first define the evolution rule for an SI model as follows:

\textbf{Definition:} For a cell $x$ in a space $\mathcal{L}$ we define the SI rule as:
\begin{equation}
    \phi_{SI}^t(x)=\left\{\begin{array}{ll}
        S & \text{if }\pi^t(x)=S\text{, }\sum_g{\sigma_{g,I}^t(x)\cdot P(g)}\leq \cdots\\
         &  \cdots\sum_g{\sigma_{g,S}^t(x)\cdot P(g)}\text{ and }\rho\geq i(t),\\
        I & \text{if }\pi^t(x)\in\{S,I\}\text{,} \\
        \pi^t(x) & \text{otherwise,}
    \end{array}\right.
\end{equation}
with $\rho\in\mathcal{U}_{[0,1]}$, $P(g)$ the impact rate of degree $g$ and 
\begin{equation}
    i(t) = \frac{\beta}{\alpha}\sum_g{\frac{\sigma_{g,I}^t(x)}{\Delta_g}}\cdot P(g),
\end{equation}
where the factor $\frac{\beta}{\alpha}$ indicates the affectation caused by the disease and the remaining factor corresponds to the behavior of the cells and its impact on the cell $x$.

The above definition allows us to naturally define the evolution rules for SIS and SIR models in cellular automata:

\textbf{Definition.} Given a cell $x$ in a set $\mathcal{L}$ we define respectively the evolution rules for the SIS and SIR models respectively as:
\begin{equation}
    \phi_{SIS}^t(x)=\left\{\begin{array}{ll}
        \phi_{SI}^t(x) & \text{if }\pi^t(x) = S,\\
        I & \text{if }\pi^t(x)=I\text{ and }\rho>\alpha,\\
        S & \text{if }\pi^t(x)=I\text{ and }\rho\leq\alpha.
    \end{array}\right.
\end{equation}

\begin{equation}
    \phi_{SIR}^t(x)=\left\{\begin{array}{ll}
        \phi_{SI}^t(x) & \text{if }\pi^t(x) = S,\\
        I & \text{if }\pi^t(x)=I\text{ and }\rho>\alpha,\\
        R & \text{if }\pi^t(x)=I\text{ and }\rho\leq\alpha, and \\
        R & \text{if }\pi^t(x)=R,
    \end{array}\right.
\end{equation}
where $\rho\in\mathcal{U}_{[0,1]}$.

For the models with birth and mortality, we will start from the hypothesis that each individual in a system has a probability of dying from non-disease causes that depends on its own age.

To implement the notion of age of a cell, we will initially have to modify the domain of our evolution rule. Specifically, we will have to take into account the state of the central cell together with its age and the state of its neighbors, so that 
\begin{equation*}
\begin{array}{c}
    Dom(\phi)=\Sigma_x\times K\times \overbrace{\Sigma\times\cdots\times\Sigma}^N, \\
    \text{with }K=\{1,2,\cdots,100\},
\end{array}
\end{equation*}
if we assume that the age of the cell $x$ can range from 1 to 100 time units (weeks, months, years, etc.).

Another feature that we will implement will be that of cell aging, this with the objective of analyzing the impact of a disease on the individuals of a system at different stages of their "life". The way we will approach this idea will be with the following adjustment to the range of our rule
$$Ran(\phi)=\Sigma_x\times K. $$
Since we are working on populations of constant size, the way we will interpret the birth of a cell will be with the occupation of the space left by one that dies. For this we will identify the spaces left by the cells that die with state $D$ and zero age, this will allow us to separate the spaces that can be occupied and those that cannot be occupied from cells that interact with their neighbors.

With these ideas in mind, we can define the evolution rule for an $M$ epidemiological model as follows:

\textbf{Definition.} Let $x$ be a cell in a set $\mathcal{L}$, $M$ be an epidemiological model ($SIS, SIR$, etc.) and $T$ be a time unit (days, months, years, etc.). We define the evolution rule with births and deaths for $M$ as:
\begin{equation}
    \mu_{M,T}^t(x)=\left\{\begin{array}{ll}
        D,0 & \text{if }t\not\equiv 0 \text{ (mod }T\text{), }\\
        & \pi^t(x)\in\{S,I,R\} \\ 
        & \text{ and }\rho\leq\omega_k, \\
        D,0 & \text{if }\pi^t(x)=D\text{ and }\rho>b,\\
        S,1 & \text{if }\pi^t(x)=D\text{ and }\rho\leq b,\\
        \phi_M^t(x),E^t(x) & \text{if }t\not\equiv 0 \text{ (mod }T),\\
        \phi_M^t(x),E^t(x)+1 & \text{if }t\equiv 0 \text{ (mod }T),
    \end{array}\right.
\end{equation}
where $\omega_k$ is the probability of dying from non-disease causes for ages in the k-th partition of the interval $[0,100]$, $b$ is the birth rate, $\phi_M^t$ is the evolution rule of the epidemiological model, $E^t(x)$ denotes the age of the cell $x$ at time $t$ and $\rho\in\mathcal{U}_{[0,1]}$.

$\textit{\underline{Observation:}}$ The probability distribution with which $\rho$ is chosen may not be uniform and instead depend on the phenomenon being modeled. For the purposes of this project, we will focus only on the uniform distribution and hope to further analyze the level of impact that the choice of different probability distributions may have in future research.

Finally, for models with death by disease we will partition over the interval $[0,100]$ and define a probability $\theta_k$ that indicates the probability of dying from the disease in the k-th partition, similar to what was done for models with birth and mortality. Thus,

\textbf{Definition.} Let $x$ be a cell in a set $\mathcal{L}$, $M$ be an epidemiological model ($SIS$, $SIR$, etc.) and $T$ be a temporal unit (days, months, years, etc.). We define the evolution rule with death by disease for $M$ as:
\begin{equation}
    \theta_{M,T}^t(x)=\left\{\begin{array}{ll}
        D,0 & \text{if }\pi^t(x)=I\text{ and }\rho\leq\theta_k, \\
        \mu_{M,T}^t(x) & \text{otherwise.}
    \end{array}\right.
\end{equation}
Where $\theta_k$ is the probability of dying from the disease for individuals with an age in the k-th interval of the interval partition $[0,100]$, $\mu_{M,T}^t$ is the evolution rule for models with births and deaths and $\rho\in\mathcal{U}_{[0,1]}$.

\section{APPLIED EXAMPLE}

Suppose we discover the beginning of a flu outbreak in the hospital of a village with 49 people, and we want to analyze the behavior of the disease in order to establish control measures for future variations of the same flu.

In the village we can identify five different places where the cells interact with each other: the school, the offices, the market, the hospital and the homes of each person. We will assume that a census was taken in the population, and we have access to the ages of each individual, their place of work (or where they spend most of their time) and the people they live with. They are distributed as follows:

\begin{itemize}
    \item \textbf{School (E):} It is known that there are 9 children and 2 teachers in the village.
    \item \textbf{Offices (O):} It has a staff of 16 individuals.
    \item \textbf{Market (M):} 8 workers were identified.
    \item \textbf{Hospital (H):} Between doctors, nurses and patients a number of 14 individuals are identified.
\end{itemize}

We will assume that all children in the village are only children and that households with three individuals have at least one child. In addition to this and according to the above information, we will assume that in the village the population lives in one of three types of housing: the C1 housing where only one adult lives; the C2 housing where two individuals live, either two adults or one adult and one child; and finally the C3 housing, where 3 people live: two adults and one child.

\begin{figure}[ht]
  \centering
    \includegraphics[width=0.3\textwidth]{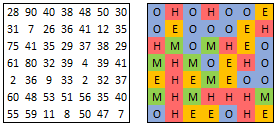}
    \caption{Ages and occupations in the village.}
    \label{figure:agesYOoccupations}
\end{figure}

\begin{figure}[ht]
  \centering
    \includegraphics[width=0.45\textwidth]{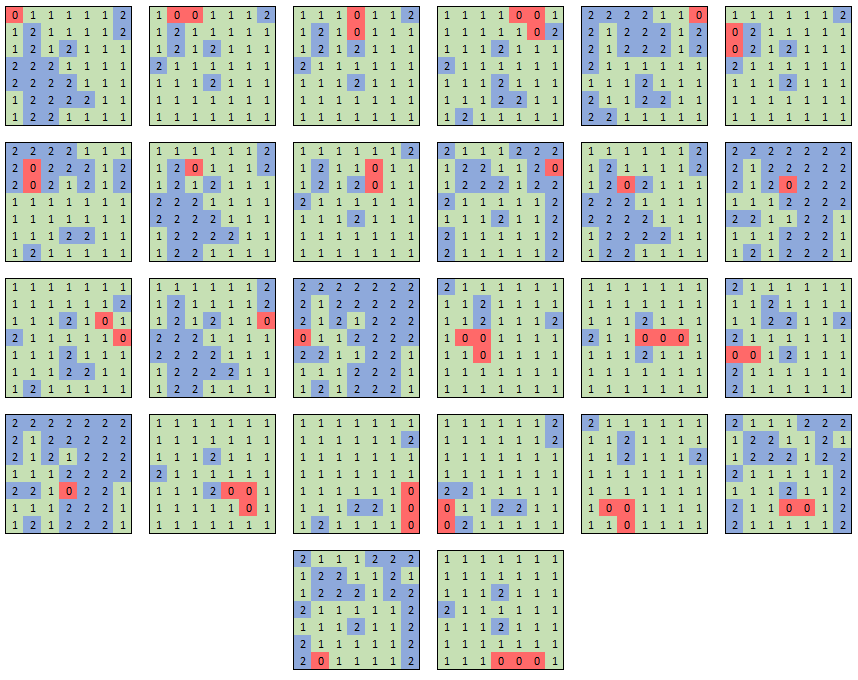}
    \caption{Impact grades by minimal neighborhood.}
    \label{fig:vecindadesEx4}
\end{figure}

Now we must establish the degrees of impact according to the way in which our population is distributed. We will assume that people who live together have an impact degree of zero among them, that with those who have contact in the place of occupation they have an impact degree equal to one and that with the other cells they have an impact degree of two. Figure \ref{fig:vecindadesEx4} describes the different neighborhood systems along with the degrees of impact generated for each cell according to the type of house it lives in (C1, C2 or C3) and the occupancy of it and its closest neighbors. 

\newpage

In the case of the disease, we will assume that it lasts an average of 5 days, so we will take $\alpha=\frac{1}{5}=0.2$ and, on the other hand, we will take an infection rate $\beta=0.3$. We will assume birth and death rates equal to 0.02 and 0.005 respectively and finally, we will assume from previous influenza outbreaks that the case fatality rates of the disease are:

\begin{table}[h]
\begin{center}
\begin{tabular}{| c | c |}
\hline
Age range & Rates \\ \hline
1 - 15 & 0.005 \\
16 - 48 & 0.01 \\
49 - 55 & 0.1 \\
56+ & 0.25 \\\hline
\end{tabular}
\caption{Estimated flu case fatality rates in village.}
\end{center}
\end{table}

We will analyze disease behavior over a 30-day period, taking different impact rates. This will give us an insight into the implications of close relationships with individuals outside the minimal neighborhood in the event of disease spread.

\begin{figure}[ht]
  \centering
    \includegraphics[width=0.45\textwidth]{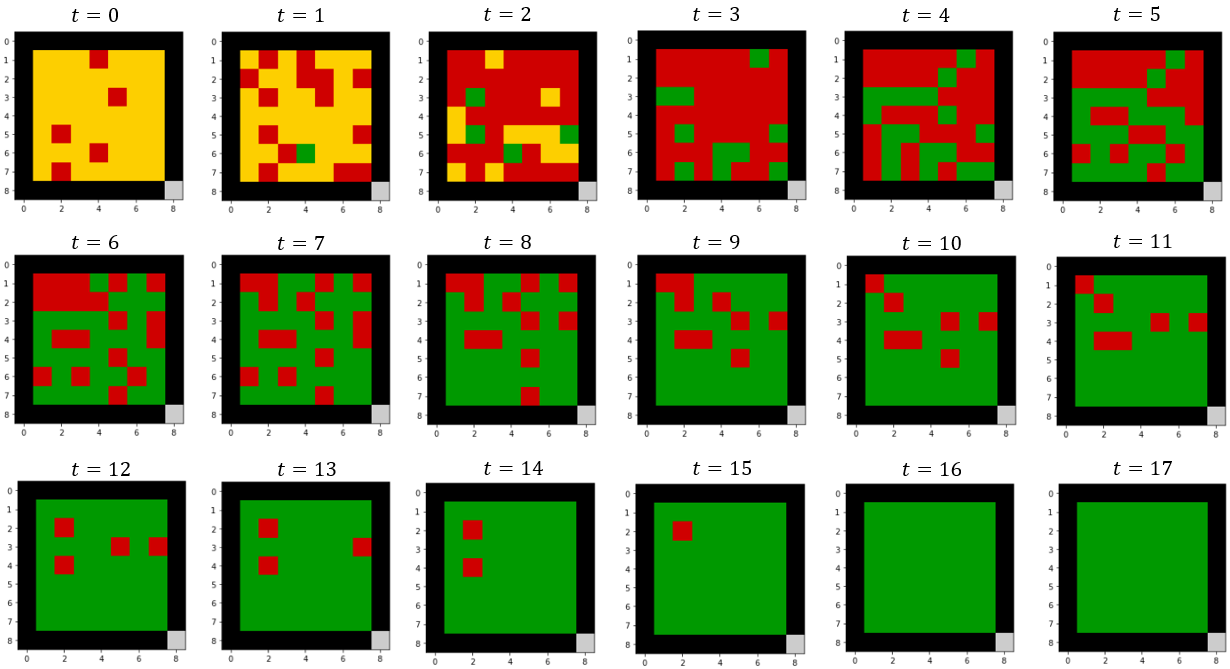}
    \caption{Evolution of the disease, taking the hospital as the starting point.}
    \label{fig:evo1}
\end{figure}

\begin{figure}[ht]
  \centering
    \includegraphics[width=0.35\textwidth]{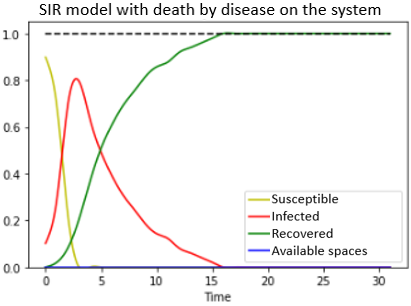}
    \caption{Changes of state with starting point in the hospital.}
    \label{fig:metricas1}
\end{figure}

In the exercise shown in figure \ref{fig:evo1} it can be observed that the disease ends up disappearing in a period of 16 days, it should be noted that at some point the entire population contracted the disease and that the number of cases began to decrease from day 3. Figure \ref{fig:metricas1} shows the number of individuals per state for a period of 30 days, taking as impact rates for grades $0,1$ and $2$ the values $P(0)=1$, $P(1)=0.5$ and $P(2)=0.25$.

If we now instead assume that within the village the relationships are a bit more distant due to e.g. measures such as alternation periods, we can assume that the impact rates are lower than in the previous exercise. In this scenario we will take $P(0)=1$, $P(1)=0.005$ and $P(2)=0.0025$ and thus,

\begin{figure}[ht]
  \centering
    \includegraphics[width=0.35\textwidth]{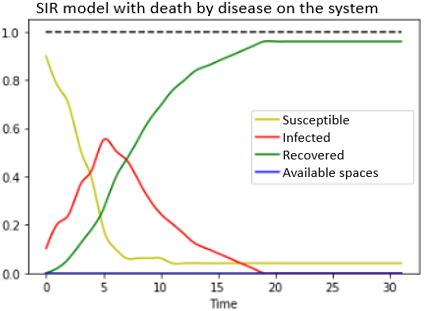}
    \caption{Changes of state with starting point in the hospital.}
    \label{fig:metricas2}
\end{figure}

It is worth noting that despite the fact that the population is more distant, the disease has a similar behavior to that of the first scenario. This is due to the nature described by the parameters of the disease and its basic reproduction number $\mathcal{R}_0\approx1.46>1$ (\ref{eq:R0}). In the figure \ref{fig:comparacionTasasdeImpacto} we can observe the average evolution over 100 simulations of the disease taking into account different impact rates:

\begin{figure}[ht]
  \centering
    \includegraphics[width=0.475\textwidth]{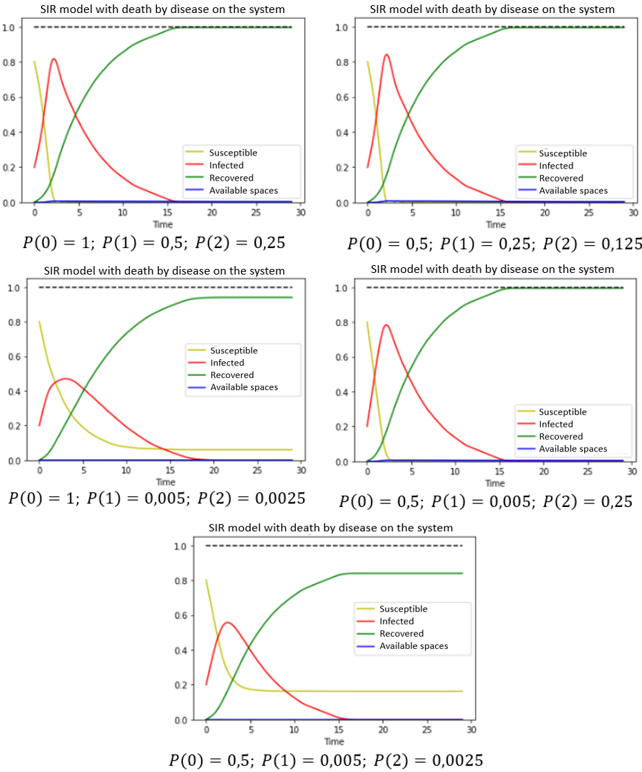}
    \caption{Average evolution of the disease taking different impact rates.}
    \label{fig:comparacionTasasdeImpacto}
\end{figure}

\begin{figure}[ht]
  \centering
    \includegraphics[width=0.4\textwidth]{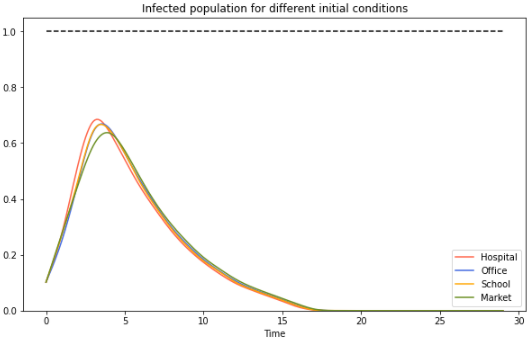}
    \caption{Infected population by initial condition.}
    \label{fig:condicionesIniciales}
\end{figure}

Considering that the previous example shows that different impact rates affect the evolution curves of the model, we can ask ourselves if something similar happens with different initial conditions. Figure \ref{fig:condicionesIniciales} shows the evolution of the infected population in the village for different initial conditions, each one located in one of the work sites mentioned above.

So far we have been working around a similar scenario with a value of $\mathcal{R}_0>1$ and in all the exercises we have been dealing with, it is evident that the disease reaches a peak (or maximum point of number of active infections) at a similar time. This is because despite considering the same parameters of the disease, we are modifying the properties of the cell space itself, and this in turn allows us to state that the initial condition will not depend only on the parameters $\alpha$, $\beta$ and the rates $\mu_k, \theta_k$ as in the classical model, but instead the ways in which the cells may interact will have to be taken into account.

\section{CAsimulations}

For the development of this work, it was necessary to design and implement a Python library that would serve as a tool to develop the analyses presented throughout this document. As a tool for the reader, different \href{https://github.com/Grupo-de-simulacion-con-automatas/Prediccion-del-comportamiento-de-una-enfermedad-simulada-en-AC-con-un-algoritmo-en-RN/tree/master/Codigo}{\underline{examples}} were developed in which the concepts defined in this document are applied. Additionally, a detailed \href{https://github.com/Grupo-de-simulacion-con-automatas/Prediccion-del-comportamiento-de-una-enfermedad-simulada-en-AC-con-un-algoritmo-en-RN/tree/master/Codigo/EpidemiologicalModels}{\underline{documentation}} on each of the functions that are part of the library is provided.

\section{CONCLUSIONS}

In this work, the dynamic nature of the models based on compartments that describe the behavior of the states that characterize a propagative phenomenon in an established set was understood. Likewise, the properties of cellular automata that allow describing large-scale spatial behaviors, from basic concepts of topology such as the fundamental systems of neighborhoods, the partial order relations that define the sets $\mathcal{V}(x)$ given any topology, among others. The understanding of these concepts allowed us to establish a relationship between the interactions of the individuals of some finite system with the fundamental systems of neighborhoods, and this in turn gave us a clear vision of the incidence of the nature of the interactions and the frequencies with which they occur in the evolution of a disease.

The conclusions that can be stated from this work are:
\begin{itemize}
    \item By making use of the properties of cellular automata to describe spatial behaviors and of fundamental neighborhood systems, it is possible to model the social relationships of a given set of individuals.
    \item The initial conditions on how cells interact do not affect the equilibrium points of the curves describing the behavior of the disease affect the equilibrium points of the curves describing the behavior of the disease. However, as observed in the examples performed, changes in the initial condition can affect the rate of propagation of the same disease.
    \item There is evidence that limiting and/or reducing the intensity of social interactions is an effective measure to reduce the number of infected individuals.
    \item Although the rules proposed allow the analysis of characteristics that were not possible with the classical models, a limitation is evident in that a maximum capacity of individuals in the system is assumed.
    \item The methodology used for the design and implementation of the evolution rules described in this work provides a clear path for the definition of rules that model the behavior of more general epidemiological models.
\end{itemize}

\addtolength{\textheight}{-12cm}   
\bibliography{BibliMSc}
\bibliographystyle{abbrv}

\end{document}